\documentclass[11pt]{article}
\usepackage[T1]{fontenc}
\usepackage{lmodern}
\usepackage{amsmath,amssymb,amsthm}
\usepackage[margin=1.05in]{geometry}

\newtheorem{theorem}{Theorem}[section]
\newtheorem{lemma}[theorem]{Lemma}
\theoremstyle{definition}
\newtheorem{definition}{Definition}[section]
\theoremstyle{plain}
\newcounter{proofstep}
\numberwithin{equation}{section}

\begin{document}

\title{Rigidity of the Peres--Vir\'ag class}
\author{Lu\'{\i}s Daniel Abreu}
\date{}
\maketitle

\begin{abstract}
Consider the hyperbolic Gaussian analytic function (GAF) 
\begin{equation*}
F_{\mathbb{D}}(z)=\sum_{n=0}^{\infty }\xi _{n}z^{n},\qquad |z|<1,
\end{equation*}%
where the $\xi _{n}$ are independent standard complex Gaussian variables.
Peres-Vir\'{a}g theorem tells that the distribution of the zeros of $F_{%
\mathbb{D}}(z)$ is provided by a Determinantal Point Process (DPP). We prove
that $F_{\mathbb{D}}(z)$ \emph{is the only GAF whose zeros are determinantal}
with a Hermitian kernel, up to multiplication by a zero-free holomorphic
function and a one-to-one holomorphic change of variable, with a \
non-degeneracy condition imposed to exclude trivial cases. The proof
identifies a geometric property of GAF's with (local) determinantal zeros:
the metric defined by the first intensity must have constant curvature $-4$.
\end{abstract}

\section{Introduction}

Consider the Gaussian power series 
\begin{equation}
F_{\mathbb{D}}(z)=\sum_{n=0}^{\infty }\xi _{n}z^{n},\qquad |z|<1,
\label{eq:hardy-gaf}
\end{equation}%
where the $\xi _{n}$ are independent standard complex Gaussian variables.
This is a Gaussian analytic function (GAF) whose covariance kernel is the
reproducing kernel of the Hardy space $H^{2}(\mathbb{D})$: 
\begin{equation*}
K_{\mathbb{D}}(z,w)=\mathbb{E}[F_{\mathbb{D}}(z)\overline{F_{\mathbb{D}}(w)}%
]=\frac{1}{1-z\overline{w}},\qquad z,w\in \mathbb{D}.
\end{equation*}

A celebrated result of Peres and Vir\'{a}g~\cite[Theorem 1]{PV} shows that
the zeros of $F_{\mathbb{D}}$ form a determinantal point process (DPP) whose
correlation kernel is the reproducing kernel of the unweighted Bergman space 
$A^{2}(\mathbb{D})$~\cite[Chapter 1]{HKZ}: 
\begin{equation}
\mathcal{K}_{\mathbb{D}}(z,w)=\frac{1}{(1-z\overline{w})^{2}}\text{,}
\label{eq:bergman-kernel}
\end{equation}%
relative to normalized area measure $d\mu =\pi ^{-1}dx\,dy$. It is the
purpose of this paper to show that (\ref{eq:hardy-gaf})\emph{\ is the only
nontrivial scalar example of a GAF whose zeros are determinantal}, beyond
the restrictions and transformations of \eqref{eq:hardy-gaf} described in
the definition below of the \emph{Peres--Vir\'{a}g class.}

To state our main definition and result, we will need some definitions.
Throughout, $\Omega $ is a nonempty connected open subset of $\mathbb{C}$, $%
d\mu =\pi ^{-1}dx\,dy$, and 
\begin{equation*}
\partial =\frac{1}{2}(\partial _{x}-i\partial _{y}),\qquad \overline{%
\partial }=\frac{1}{2}(\partial _{x}+i\partial _{y}).
\end{equation*}%
A Gaussian analytic function (GAF) is a random holomorphic function $F$
whose finite collections of values are centered complex Gaussian vectors and
satisfy 
\begin{equation*}
\mathbb{E}[F(z)F(w)]=0.
\end{equation*}%
Its covariance kernel is denoted by 
\begin{equation}
K(z,w)=\mathbb{E}[F(z)\overline{F(w)}].  \label{covarianceKernel}
\end{equation}%
Thus $K$ is holomorphic in $z$ and antiholomorphic in $w$. Multiplication of %
\eqref{eq:hardy-gaf} by a zero-free holomorphic function preserves its zeros~%
\cite[proof of Theorem 2.5.2]{HKPV}, while a one-to-one holomorphic change
of variable $\phi :\Omega \longrightarrow \mathbb{D}$ preserves the
determinantal zero property in the domain $\Omega $ (see~\cite[Section 4.2]%
{HKPV} and  Section~\ref{sec:preliminaries} below). For the classical
extension to bounded simply connected domains with smooth boundary see~\cite[%
Corollary 5.4.1]{HKPV}.

\begin{definition}
\label{def:pv-class} Let $\Omega \subseteq \mathbb{C}$ be a nonempty
connected open set, possibly $\Omega =\mathbb{C}$. A Gaussian analytic
function $F$ on $\Omega $ belongs to the \emph{Peres--Vir\'{a}g class} if 
\begin{equation}
\mathrm{Cov}\bigl(F(z),F^{\prime }(z),F^{\prime \prime }(z)\bigr)>0\qquad 
\text{for every }z\in \Omega ,  \label{eq:derivative-assumption}
\end{equation}%
and its zero process is a Hermitian DPP relative to $d\mu $.
\end{definition}

The inequality in \eqref{eq:derivative-assumption} means that the covariance
matrix is positive definite. It implies that $K(z,z)>0$ and the first zero
intensity is positive. This condition excludes the trivial case of a
Gaussian linear polynomial, which has one random zero and gives a DPP (see
Section~\ref{subsec:linear-example}). Under \eqref{eq:derivative-assumption}%
, the zeros are almost surely isolated and simple~\cite[Lemma 2.4.1]{HKPV}.

\begin{theorem}
\label{thm:main} A Gaussian analytic function $F$ on $\Omega $ belongs to
the Peres--Vir\'{a}g class if and only if there are a zero-free holomorphic
function $h$ on $\Omega $ and a one-to-one holomorphic map $\phi :\Omega
\longrightarrow \mathbb{D}$ such that 
\begin{equation}
K(z,w)=\frac{h(z)\overline{h(w)}}{1-\phi (z)\overline{\phi (w)}}.
\label{eq:main-covariance}
\end{equation}%
In this case, the zero process admits a correlation kernel relative to $d\mu 
$, holomorphic in the first variable and antiholomorphic in the second,
given by 
\begin{equation}
\mathcal{K}_{\phi }(z,w)=\frac{\phi ^{\prime }(z)\overline{\phi ^{\prime }(w)%
}}{(1-\phi (z)\overline{\phi (w)})^{2}}.  \label{eq:main-dpp-kernel}
\end{equation}
\end{theorem}

The proof is organized along Section~\ref{sec:proof} and will be outlined
after the following statement of Liouville's theorem, which can be found in~%
\cite[Theorem C]{KR} or in~\cite[Theorem 9.8]{BM}.

\begin{lemma}[Liouville's theorem]
\label{lem:liouville} Let $U\subset \mathbb{C}$ be simply connected, and let 
$g>0$ be smooth with 
\begin{equation}
\partial \overline{\partial }\log g(z)=2g(z).  \label{Curvature}
\end{equation}%
There is a holomorphic map $\phi :U\rightarrow \mathbb{D}$ with $\phi
^{\prime }(z)\neq 0$ such that 
\begin{equation}
g(z)=\frac{|\phi ^{\prime }(z)|^{2}}{(1-|\phi (z)|^{2})^{2}}.
\label{eq:liouville-map}
\end{equation}
\end{lemma}

The proof of Theorem \ref{thm:main} relies on a clear mechanism that we
proceed to describe. Edelman--Kostlan's formula~\cite[Section 8]{EK} for the 
$1$-point zero intensity of a GAF with covariance kernel $K(z,w)$: 
\begin{equation}
\rho _{1}(z)=\partial \overline{\partial }\log K(z,z).  \label{eq:gQ}
\end{equation}

First recall that Sodin's Calabi rigidity theorem~\cite[Theorem 2]{Sodin},%
\cite[Theorem 2.5.2]{HKPV} determines the GAF from $\rho _{1}(z)$, up to
multiplication by a zero-free holomorphic function. Consequently, if we show
that $\rho _{1}(z)\ $satisfies (\ref{Curvature}) in $\Omega $ (not assumed
to be simply connected), Liouville's theorem gives (\ref{eq:liouville-map}),
which is is the diagonal of the kernel (\ref{eq:main-dpp-kernel}), proving a
restriction of Theorem \ref{thm:main}, on an open disc $U\subset \Omega $,
with a holomorphic map $\phi :U\rightarrow \mathbb{D}$ with $\phi ^{\prime
}(z)\neq 0$. Somehow surprisingly, the identities obtained by imposing the
determinantal structure on diagonal perturbations of the $2$- and $3$-point
Kac-Rice intensities are enough to prove (\ref{Curvature}), after showing
that they impose a holomorphic structure on the determinantal kernel, which
is first proved to be smooth. The local version of Theorem \ref{thm:main} is
finally extended to $\Omega $ by analytic continuation and $\phi $ is shown
to be globally one-to-one (this only uses the global repulsiveness of the $2$%
-point intensity).

Our results therefore identify the local geometric mechanism for membership
in the Peres-Vir\'{a}g class: the local determinantal assumption forces the
metric $ds^{2}=\partial \overline{\partial }\log K(z,z)\left\vert
dz\right\vert ^{2}$ to have constant curvature $-4$ (this is equivalent to (%
\ref{Curvature})). Moreover, Definition \ref{def:pv-class} can be relaxed to
hold locally, provided the GAF zeros have a global two-point intensity such
that (\ref{eq:G}) satisfies $G(z,w)\leq 1$ in the whole domain. In Section
4.4 we construct examples of GAF's satisfying the curvature condition, with
determinantal zeros locally, but not globally. This shows that global
injectivity in Theorem~\ref{thm:main} cannot be replaced by local
injectivity. A limit case of Katori--Shirai's GAF~\cite{KS} is globally
repulsive but fails the curvature condition (see Section 4.3).

There is a considerable body of work motivated by the work of Peres and Vir%
\'{a}g. We list a few developments to settings beyond the scope of our
results. A real-coefficient counterpart is provided by the limiting Kac
series, with independent standard \emph{real} Gaussian coefficients, whose
zero correlations admit Pfaffian formulas. Forrester~\cite{Forrester}
obtained these through truncated random orthogonal matrices, while Matsumoto
and Shirai~\cite[Theorems 2.1 and 2.7]{MS} gave a direct proof for real
zeros on $(-1,1)$ and nonreal zeros in the upper half-disk. Matrix
constructions with determinantal zero processes have been developed by
Krishnapur~\cite[Theorem 4]{Krishnapur}.

The outline of the paper is straightforward. The preliminary Section~\ref%
{sec:preliminaries} recalls the facts about DPPs and GAFs used in the proof.
The proof of Theorem~\ref{thm:main} is given in Section~\ref{sec:proof}: two
preparation lemmas are followed by the curvature lemma and a two-step proof
of the theorem. We conclude with examples and comments in Section~\ref%
{sec:examples}.

\section{Preliminaries on DPPs and GAFs}

\label{sec:preliminaries}

\subsection{Hermitian determinantal point processes}

Let $\mathcal{Z}$ be a random subset of $\Omega $, with no repeated points
and finitely many points in every compact subset. Its $n$-point intensity $%
\rho _{n}$, relative to the normalized are measure $d\mu =\pi ^{-1}dx\,dy$,
is defined by 
\begin{equation*}
\mathbb{E}\sum_{z_{1},\ldots ,z_{n}\in \mathcal{Z}}^{\neq }f(z_{1},\ldots
,z_{n})=\int_{\Omega ^{n}}f(z_{1},\ldots ,z_{n})\rho _{n}(z_{1},\ldots
,z_{n})\,d\mu (z_{1})\cdots d\mu (z_{n})
\end{equation*}%
for nonnegative measurable $f$. The sum is over ordered tuples of distinct
points. Integrating $\rho _{1}$ over a set gives its expected number of
points, and integrating $\rho _{2}$ over a subset of $\Omega \times \Omega $
gives its expected number of ordered pairs of distinct points.

For a Hermitian DPP we use the usual operator definition. The correlation
operator $\mathcal{L}$ on $L^{2}(\Omega ,d\mu )$ is selfadjoint and
satisfies $0\leq \mathcal{L}\leq I$. For every relatively compact set $U$,
the operator $\mathbf{1}_{U}\mathcal{L}\mathbf{1}_{U}$ has finite trace.
Here $\mathbf{1}_{U}$ denotes multiplication by the indicator of $U$. We
denote the measurable Hermitian integral kernel of $\mathcal{L}$ by $%
\mathcal{K}(z,w)$, with the diagonal interpreted through the local spectral
expansion used in the proof of Lemma~\ref{lem:local-factorization}. These
are the standard operator assumptions for Hermitian DPPs; see~\cite%
{HKPV,Borodin,Soshnikov}. The intensities are then given by 
\begin{equation}
\rho _{n}(z_{1},\ldots ,z_{n})=\det [\mathcal{K}(z_{i},z_{j})]_{i,j=1}^{n}%
\quad \text{almost everywhere}.  \label{eq:dpp}
\end{equation}%
The determinant for two points immediately gives 
\begin{equation}
|\mathcal{K}(z,w)|^{2}=\rho _{1}(z)\rho _{1}(w)-\rho _{2}(z,w).
\label{eq:modulus-identity}
\end{equation}%
Multiplying the kernel by the factors 
\begin{equation*}
\mathcal{K}(z,w)\longmapsto e^{i\theta (z)}\mathcal{K}(z,w)e^{-i\theta (w)}
\end{equation*}%
and taking the complex conjugate of the kernel, leaves every determinant
unchanged. Restriction to a domain $U\subseteq \Omega $ preserves the DPP
property~\cite[Remark 4.2.5]{HKPV}. If $\psi :V\longrightarrow U$ is a
one-to-one holomorphic map from a planar domain $V$ onto $U$, a kernel
relative to $d\mu $ for the process $\psi ^{-1}(\mathcal{Z}\cap U)$ on $V$
is 
\begin{equation*}
\psi ^{\prime }(z)\overline{\psi ^{\prime }(w)}\mathcal{K}(\psi (z),\psi
(w)),\qquad z,w\in V.
\end{equation*}%
Consequently the normalized pair correlation function, defined as the ratio 
\begin{equation}
G(z,w)=\frac{\rho _{2}(z,w)}{\rho _{1}(z)\rho _{1}(w)},  \label{eq:G}
\end{equation}%
is the same in either coordinate, after composition with the coordinate map.

\subsection{Gaussian analytic functions and Kac--Rice intensities}

A standard complex Gaussian variable $\xi $ has density $e^{-|z|^{2}}$
relative to $d\mu $ on $\mathbb{C}$. Its real and imaginary parts are
independent, each with variance $1/2$.  Derivatives of a GAF are again
jointly Gaussian with its values. In \eqref{eq:gQ}, the derivatives of $%
K(z,z)$ are taken after restriction to the diagonal. For any jointly
centered circular complex Gaussian vectors $Y,\xi $, with $\mathrm{Cov}(\xi )
$ invertible, the Gaussian conditioning identity is 
\begin{equation}
\mathrm{Cov}(Y\mid \xi =0)=\mathrm{Cov}(Y)-\mathrm{Cov}(Y,\xi )\mathrm{Cov}%
(\xi )^{-1}\mathrm{Cov}(\xi ,Y).  \label{eq:gaussian-conditioning}
\end{equation}%
The first intensity is provided by Edelman-Kostlan formula (\ref{eq:gQ}).\
We will follow~\cite[Section 3.4, Corollary 3.4.2]{HKPV} in the presentation
of the Kac-Rice formulas providing the higher intensities. For distinct
points $z_{1},\ldots ,z_{n}$ at which the covariance of the values is
nonsingular, put 
\begin{equation*}
A_{ij}=K(z_{i},z_{j}),\quad B_{ij}=\mathbb{E}[F^{\prime }(z_{i})\overline{%
F(z_{j})}],\quad C_{ij}=\mathbb{E}[F^{\prime }(z_{i})\overline{F^{\prime
}(z_{j})}].
\end{equation*}%
By \eqref{eq:gaussian-conditioning}, conditioning on the values being zero
gives a centered Gaussian vector of derivatives, whose covariance, 
\begin{equation}
\Lambda =\mathrm{Cov}\!\left( (F^{\prime }(z_{i}))_{i=1}^{n}\,\mid
\,F(z_{1})=\cdots =F(z_{n})=0\right) =C-BA^{-1}B^{\ast },
\label{eq:conditional-covariance}
\end{equation}%
is called the \emph{conditional covariance matrix}. The Kac--Rice formula
providing the $n$-point intensities of the zeros of a GAF is 
\begin{equation}
\rho _{n}(z_{1},\ldots ,z_{n})=\frac{\mathbb{E}[\prod_{j=1}^{n}|F^{\prime
}(z_{j})|^{2}\mid F(z_{1})=\cdots =F(z_{n})=0]}{\det A}=\frac{\mathrm{per}%
\,\Lambda }{\det A}.  \label{eq:kac-rice}
\end{equation}%
Here $\mathrm{per}\,M=\sum_{\sigma }\prod_{i}M_{i,\sigma (i)}$, where the
sum is over all permutations of $\{1,\ldots ,n\}$. The second equality
follows by computing the Gaussian moment in the numerator~\cite[Proposition 8%
]{PV}.

For a GAF in the Peres--Vir\'{a}g class, the determinantal assumption
implies that the same intensity has two expressions: 
\begin{equation}
\rho _{n}(z_{1},\ldots ,z_{n})=\frac{\mathrm{per}\,\Lambda }{\det A}=\det [%
\mathcal{K}(z_{i},z_{j})]_{i,j=1}^{n},  \label{eq:joint-intensity-identities}
\end{equation}%
where Kac--Rice applies, with the determinantal equality initially holding
almost everywhere.

\section{Proof of Theorem~\protect\ref{thm:main}}

\label{sec:proof}

\subsection{Perturbation of diagonal Kac--Rice intensities}

\label{sec:preparation}

To study the intensities near the diagonal, we need some identities
involving the covariance matrices of the values and derivatives and their
determinants, evaluated at the same point $z$: 
\begin{equation*}
J_{m}(z)=[\mathbb{E}(F^{(i)}(z)\overline{F^{(j)}(z)})]_{i,j=0}^{m},\qquad
D_{m}(z)=\det J_{m}(z).
\end{equation*}%
The required determinant identities needed below are 
\begin{eqnarray}
D_0(z) &=& K(z,z),\qquad D_1(z)=K(z,z)^2\rho_1(z),  \notag \\
D_2(z) &=& K(z,z)^3\rho_1(z)^2 \bigl(2\rho_1(z)+\partial\overline{\partial}%
\log\rho_1(z)\bigr).  \label{eq:derivative-determinants}
\end{eqnarray}%
For the last identity in \eqref{eq:derivative-determinants} combine $\log
D_{1}(z)=2\log K(z,z)+\log \rho _{1}(z)$ with Jacobi's determinant identity, 
\begin{equation*}
\partial \overline{\partial }\log D_{1}(z)=\frac{D_{2}(z)D_{0}(z)}{%
D_{1}(z)^{2}}.
\end{equation*}%
The conditional variance of $F^{\prime \prime }(z)$ is obtained by
conditioning on the two Gaussian variables $F(z),F^{\prime }(z)$ at the same
point. Write 
\begin{eqnarray*}
&&J_2(z)=\left(%
\begin{array}{cc}
J_1(z) & b(z)^* \\ 
b(z) & c(z)%
\end{array}%
\right),\qquad c(z)=\mathbb{E}|F^{\prime \prime }(z)|^2, \\
&&b(z)=\bigl(\mathbb{E}[F^{\prime \prime }(z)\overline{F(z)}], \mathbb{E}%
[F^{\prime \prime }(z)\overline{F^{\prime }(z)}]\bigr).
\end{eqnarray*}%
The Gaussian conditioning formula \eqref{eq:gaussian-conditioning} gives the
variance $c(z)-b(z)J_{1}(z)^{-1}b(z)^{\ast }$, while the block determinant
formula gives 
\begin{equation*}
D_{2}(z)=D_{1}(z)\bigl(c(z)-b(z)J_{1}(z)^{-1}b(z)^{\ast }\bigr).
\end{equation*}%
Together with \eqref{eq:derivative-determinants}, this yields 
\begin{eqnarray}
V_2(z) &:=& \mathrm{Var}(F^{\prime \prime }(z)\mid F(z)=F^{\prime }(z)=0)=%
\frac{D_2(z)}{D_1(z)}  \notag \\
&=& K(z,z)\rho_1(z)\bigl(2\rho_1(z)+\partial\overline{\partial}\log \rho_1(z)%
\bigr).  \label{eq:second-conditional}
\end{eqnarray}%
Assumption \eqref{eq:derivative-assumption} makes both $\rho _{1}(z)$ and $%
2\rho _{1}(z)+\partial \overline{\partial }\log \rho _{1}(z)$ strictly
positive.

\begin{lemma}
\label{lem:nearby-intensities}\label{step:nearby} Let $F$ be a GAF on $%
\Omega $ satisfying \eqref{eq:derivative-assumption}. Then $\rho_2$ extends
real-analytically across $z=w$, with 
\begin{equation}  \label{eq:two-nearby}
\rho_2(z,z+\delta) =\frac{\rho_1(z)(2\rho_1(z)+\partial\overline{\partial}%
\log \rho_1(z))^2} {8}|\delta|^2+O(|\delta|^3).
\end{equation}
Moreover, $\rho_3$ is continuous at nearby distinct triples and, for fixed
distinct $u_1,u_2,u_3$, 
\begin{equation}  \label{eq:three-nearby}
\rho_3(z+\epsilon u_1,z+\epsilon u_2,z+\epsilon u_3) =O(|\epsilon|^6).
\end{equation}
\end{lemma}

\textbf{Proof. }The mechanism of the proof is inspired by the near-diagonal
Kac--Rice argument of Peres-Vir\'{a}g~\cite[Lemma 32]{PV} and Michelen-Yakir 
\cite[Section 3.2]{MY}.  We use divided differences for both estimates,
replacing nearly identical values by quantities whose covariance remains
nonsingular as the points approach one another. Conditioning on the Kac-Rice
conditions $F(z_{1})=\cdots =F(z_{n})=0$, we express each derivative as a
divided difference multiplied by the product of the corresponding
differences between the points (this product also appears in the two-sided
bounds of \cite[Theorem 1.1]{NS}, under stronger nondegeneracy assumptions).
Let $F[z_{1}]=F(z_{1})$ and 
\begin{equation*}
F[z_{1},\ldots ,z_{m}]=\frac{F[z_{2},\ldots ,z_{m}]-F[z_{1},\ldots ,z_{m-1}]%
}{z_{m}-z_{1}}.
\end{equation*}%
When nodes coincide, the value is the analytic limit. In particular, 
\begin{equation*}
F[z,w]=\frac{F(w)-F(z)}{w-z},\qquad F[z,z]=F^{\prime }(z),\qquad F[z,z,z]=%
\frac{F^{\prime \prime }(z)}{2}.
\end{equation*}%
Now replace 
\begin{equation*}
(F(z),F(w))\quad \text{by}\quad (F[z],F[z,w]).
\end{equation*}%
The new pair tends to $(F(z),F^{\prime }(z))$ as $w\rightarrow z$, whose
covariance is nonsingular by \eqref{eq:derivative-assumption}. For $w\neq z$%
, the conditions $F(z)=F(w)=0$ and $F[z]=F[z,w]=0$ are equivalent. The
linear change of variables has determinant $(w-z)^{-1}$, so the Kac--Rice
formula \eqref{eq:kac-rice} gives 
\begin{equation}
\rho _{2}(z,w)=\frac{\mathbb{E}[|F^{\prime }(z)F^{\prime }(w)|^{2}\mid
F[z]=F[z,w]=0]}{|w-z|^{2}\det \mathrm{Cov}(F[z],F[z,w])}.
\label{eq:pair-differences}
\end{equation}%
The divided-difference identities give 
\begin{eqnarray}
F^{\prime }(z) &=&F[z,w]-(w-z)F[z,z,w],  \notag \\
F^{\prime }(w) &=&F[z,w]+(w-z)F[z,w,w].
\label{eq:divided-difference-identities}
\end{eqnarray}%
Under the conditioning $F[z]=F[z,w]=0$, the scaled derivatives therefore
satisfy 
\begin{equation}
\left( \frac{F^{\prime }(z)}{w-z},\frac{F^{\prime }(w)}{w-z}\right) =\bigl(%
-F[z,z,w],F[z,w,w]\bigr).  \label{scaled}
\end{equation}%
Writing $\delta =w-z$, Taylor's formula gives 
\begin{eqnarray*}
F[z,w] &=&F^{\prime }(z)+\frac{\delta }{2}F^{\prime \prime }(z)+O(|\delta
|^{2}), \\
F^{\prime }(w) &=&F^{\prime }(z)+\delta F^{\prime \prime }(z)+O(|\delta
|^{2}).
\end{eqnarray*}%
We next justify the regularity of their conditional moments at the diagonal.
Joint covariances of divided differences are obtained by taking divided
differences of $K$ in its first variable and conjugate divided differences
in its second. For example, 
\begin{equation*}
\mathbb{E}[F[z]\overline{F[z,w]}]=\frac{K(z,w)-K(z,z)}{\overline{w}-%
\overline{z}}.
\end{equation*}%
Since $K$ is holomorphic--antiholomorphic, the apparent singularity at $w=z$
is removable. The same argument applies to all the entries of the covariance
of $(F[z],F[z,w],F[z,z,w],F[z,w,w])$, which extend jointly as real-analytic
functions of $z,w$ near the diagonal. Apply the Gaussian conditioning
identity \eqref{eq:gaussian-conditioning} with 
\begin{eqnarray*}
A(z,w) &=&\mathrm{Cov}(F[z],F[z,w]), \\
C(z,w) &=&\mathrm{Cov}(F[z,z,w],F[z,w,w]),
\end{eqnarray*}%
and with $B(z,w)$ the cross-covariance of the second pair with the first. On
the diagonal $A(z,z)=J_{1}(z)>0$. Thus $\det A(z,w)$ stays positive nearby,
and the cofactor formula shows that $A(z,w)^{-1}$ is real-analytic there. As
a result, the conditional covariance matrix of $\left(
F[z,z,w],F[z,w,w]\right) $ is 
\begin{equation*}
\Lambda (z,w)=C(z,w)-B(z,w)A(z,w)^{-1}B(z,w)^{\ast },
\end{equation*}%
with real-analytic entries. The corresponding conditional fourth moment is 
\begin{eqnarray*}
&&\mathbb{E}[|F[z,z,w]F[z,w,w]|^{2}\mid F[z]=F[z,w]=0] \\
&&\qquad =\Lambda _{11}(z,w)\Lambda _{22}(z,w)+|\Lambda _{12}(z,w)|^{2}.
\end{eqnarray*}%
Since this is a polynomial in the covariance entries, it is real-analytic as
well. Each conditioned derivative contributes a factor $w-z$. Consequently,
the product of their squared moduli contributes the factor $|w-z|^{4}$, and
the numerator in \eqref{eq:pair-differences} is $|w-z|^{4}$ times a
real-analytic function. The denominator contains $|w-z|^{2}$ times a
real-analytic covariance determinant whose value at $w=z$ is $D_{1}(z)>0$.
Therefore $\rho _{2}(z,w)/|w-z|^{2}$ extends real-analytically across the
diagonal. We now compute the leading quadratic coefficient of $\rho _{2}(z,w)
$ (its exact value will be later used to identify $p(z)$ in %
\eqref{eq:compatibility}). At $w=z$, the conditioning pair becomes $%
(F(z),F^{\prime }(z))$. Substitution into the divided-difference identities %
\eqref{eq:divided-difference-identities} shows that both $F[z,z,w]$ and $%
F[z,w,w]$ tend to $F^{\prime \prime }(z)/2$. By \eqref{scaled}, the two
scaled derivatives have limits $-F^{\prime \prime }(z)/2$ and $F^{\prime
\prime }(z)/2$, respectively. The product of the squared moduli of the two
limiting scaled derivatives is $\frac{|F^{\prime \prime }(z)|^{4}}{16}$.
Thus the value of the remaining real-analytic function on the diagonal is 
\begin{equation*}
\frac{1}{16}\mathbb{E}[|F^{\prime \prime }(z)|^{4}\mid F(z)=F^{\prime }(z)=0]%
\text{.}
\end{equation*}%
Under this conditioning, $F^{\prime \prime }(z)$ is a centered circular
complex Gaussian variable with variance $V_{2}(z)$, defined in %
\eqref{eq:second-conditional}. Its fourth moment is twice the square of its
variance. Hence the last expression equals 
\begin{equation*}
\frac{1}{16}\,2V_{2}(z)^{2}=\frac{V_{2}(z)^{2}}{8}.
\end{equation*}%
Dividing by the limiting covariance determinant $D_{1}(z)$ gives the
diagonal value $V_{2}(z)^{2}/(8D_{1}(z))$. As two zeros approach one
another, their vanishing conditions approach $F(z)=F^{\prime }(z)=0$.
Accordingly, 
\begin{equation*}
\rho _{2}(z,z+\delta )=\frac{V_{2}(z)^{2}}{8D_{1}(z)}|\delta |^{2}+O(|\delta
|^{3}).
\end{equation*}%
Substituting \eqref{eq:derivative-determinants} and %
\eqref{eq:second-conditional} proves \eqref{eq:two-nearby}, with a strictly
positive leading coefficient.

It remains to prove \eqref{eq:three-nearby}. We again use divided
differences to keep the value covariance nonsingular and isolate the factors
arising from the distances between the points. Here only the order of
vanishing is needed, so it suffices to bound the remaining conditional
moment without computing its value. Let $\Delta =\prod_{i<j}(z_{j}-z_{i})$
and set 
\begin{equation*}
\xi =(F[z_{1}],F[z_{1},z_{2}],F[z_{1},z_{2},z_{3}]).
\end{equation*}%
The transformation from the three values to $\xi $ has determinant $\Delta
^{-1}$, so 
\begin{equation*}
\det \mathrm{Cov}(F(z_{1}),F(z_{2}),F(z_{3}))=|\Delta |^{2}\det \mathrm{Cov}%
(\xi ).
\end{equation*}%
As $z_{1},z_{2},z_{3}\rightarrow z$, the covariance of $\xi $ tends to that
of $(F(z),F^{\prime }(z),F^{\prime \prime }(z)/2)$, whose determinant is $%
D_{2}(z)/4>0$. After restricting the points to a sufficiently small disk,
this covariance remains positive definite and its inverse is bounded. If the
three values are conditioned to be zero, the interpolation identity in~\cite[%
Lemma 3.5]{MY} gives 
\begin{equation*}
F^{\prime }(z_{i})=\prod_{j\neq
i}(z_{i}-z_{j})\,F[z_{1},z_{2},z_{3},z_{i}],\qquad i=1,2,3.
\end{equation*}%
The covariance entries of these last divided differences extend
real-analytically when the nodes coincide. Applying the Gaussian conditional
covariance formula (\ref{eq:gaussian-conditioning}) to the vector 
\begin{equation*}
\left( F[z_{1},z_{2},z_{3},z_{i}]\right) _{i=1}^{3},
\end{equation*}%
with conditioning vector $\xi $, therefore gives bounded conditional
covariances for these divided differences. Since 
\begin{equation*}
\prod_{i=1}^{3}\prod_{j\neq i}|z_{i}-z_{j}|^{2}=|\Delta |^{4},
\end{equation*}%
the numerator in Kac--Rice contains $|\Delta |^{4}$ times a bounded
continuous function. Its denominator contains $|\Delta |^{2}$ times a
positive continuous function. Thus $\rho _{3}$ is $|\Delta |^{2}$ times a
bounded continuous function. In particular it is continuous for nearby
distinct triples. For $z_{i}=z+\epsilon u_{i}$, with the $u_{i}$ fixed and
distinct, 
\begin{equation*}
|\Delta |^{2}=|\epsilon |^{6}\prod_{i<j}|u_{j}-u_{i}|^{2},
\end{equation*}%
which proves \eqref{eq:three-nearby}.  \hfill $\square $

\subsection{A holomorphic factorization for the kernel of the Peres--Vir\'{a}%
g class}

\label{sec:determinantal-kernel}

The identities \eqref{eq:joint-intensity-identities} already involve a
determinantal kernel $\mathcal{K}$, whose existence is part of the Peres--Vir%
\'{a}g assumption. At this stage, however, we only know that it is
measurable. The next result shows that the process admits a correlation
kernel with the smoothness and analytic structure required for the proof of
Theorem~\ref{thm:main}. This structure follows by expanding the $2$- and $3$%
-point determinantal intensities near the diagonal and comparing these
expansions with the Kac--Rice intensities obtained in the previous section.

\begin{lemma}
\label{lem:local-factorization} Let $F$ belong to the Peres--Vir\'ag class.
On a sufficiently small disk around each point of $\Omega$, its zero process
admits a correlation kernel of the form 
\begin{equation}  \label{eq:kernel-factor}
\mathcal{K}(z,w)=a(z)\overline{a(w)}H(z,w),
\end{equation}
where $a$ is smooth and never zero, and $H$ is a positive kernel holomorphic
in its first variable and antiholomorphic in its second.
\end{lemma}

\setcounter{proofstep}{0} \textbf{Proof. }We first obtain a smooth
correlation kernel in Step 1, and then compare its two- and three-point
determinants with the estimates of Lemma~\ref{lem:nearby-intensities} in
Step 2, leading to (\ref{eq:kernel-factor}).

\medskip \refstepcounter{proofstep}\label{step:smooth} \noindent \textbf{%
Step \theproofstep.}

Recall that the DPP property is preserved by restriction to a domain \cite[%
Remark 4.2.5]{HKPV}. The local trace-class assumption \cite[Assumption 4.2.3]%
{HKPV}, gives 
\begin{equation*}
\operatorname{tr}(\mathbf{1}_{U}\mathcal{L}\mathbf{1}_{U})<\infty \qquad \text{for
every relatively compact }U\subset \Omega .
\end{equation*}%
Fix a point $z_{0}$ and take a sufficiently small relatively compact disk $U$
around it. We can choose orthonormal eigenfunctions $e_{j}$ with eigenvalues 
$\lambda _{j}\geq 0$ for the positive trace-class operator $\mathbf{1}_{U}%
\mathcal{L}\mathbf{1}_{U}$ on $L^{2}(U,d\mu )$. The kernel and its diagonal
are given, almost everywhere, by the expansions 
\begin{equation*}
\mathcal{K}(z,w)=\sum_{j}\lambda _{j}e_{j}(z)\overline{e_{j}(w)},\qquad \rho
_{1}(z)=\sum_{j}\lambda _{j}|e_{j}(z)|^{2}.
\end{equation*}%
Since $\sum_{j}\lambda _{j}<\infty $, the vector 
\begin{equation}
v(z)=(\sqrt{\lambda _{j}}\,e_{j}(z))_{j\geq 1}  \label{eq:kernel-vectors}
\end{equation}%
belongs to $\ell ^{2}$ for almost every $z$. Thus 
\begin{equation*}
\mathcal{K}(z,w)=\langle v(z),v(w)\rangle _{\ell ^{2}},\qquad \Vert
v(z)\Vert _{\ell ^{2}}^{2}=\rho _{1}(z)\text{,}
\end{equation*}%
with the inner products linear in the first variable. Equation~%
\eqref{eq:modulus-identity} gives 
\begin{equation*}
|\mathcal{K}(z,w)|^{2}=\rho _{1}(z)\rho _{1}(w)-\rho _{2}(z,w)
\end{equation*}%
almost everywhere. By the Edelman--Kostlan formula and Lemma~\ref%
{lem:nearby-intensities}, the right-hand side is smooth for every pair in a
sufficiently small disk, including the diagonal. Now let 
\begin{equation*}
X(z)=v(z)v(z)^{\ast },
\end{equation*}%
where $v(z)v(z)^{\ast }$ is the operator $x\mapsto \langle x,v(z)\rangle
_{\ell ^{2}}v(z)$. In the real Hilbert space of selfadjoint Hilbert--Schmidt
operators, whose norm satisfies $\Vert X\Vert _{\mathrm{HS}}^{2}=\mathrm{tr}%
(X^{2})$, we have 
\begin{eqnarray*}
\Vert X(z)\Vert _{\mathrm{HS}} &=&\rho _{1}(z), \\
\langle X(z),X(w)\rangle _{\mathrm{HS}} &=&\mathrm{tr}(X(z)X(w))=\rho
_{1}(z)\rho _{1}(w)-\rho _{2}(z,w)
\end{eqnarray*}%
almost everywhere. Thus $X$ is locally square-integrable and its inner
products have a smooth representative. The smoothness criterion for
Hilbert-space-valued maps therefore gives a smooth function equal to $X$
pointwise \cite[Theorem 7.1]{Neeb}. To obtain an almost-everywhere identity,
set 
\begin{equation*}
M(z,w)=\langle X(z),X(w)\rangle _{\mathrm{HS}}\quad \text{a.e.},\qquad
X_{\varepsilon }=\eta _{\varepsilon }\ast X,
\end{equation*}%
where $M$ denotes the smooth representative of the scalar kernel and $\eta
_{\varepsilon }$ is a standard real mollifier with unit integral. For
real-coordinate multi-indices $\alpha ,\beta $, Fubini's theorem and
integration by parts give 
\begin{equation*}
\big\langle D^{\alpha }X_{\varepsilon }(z),D^{\beta }X_{\delta }(w)%
\big\rangle_{\mathrm{HS}}=\iint \eta _{\varepsilon }(z-u)\eta _{\delta
}(w-v)D_{u}^{\alpha }D_{v}^{\beta }M(u,v)\,du\,dv.
\end{equation*}%
Here $z,w$ lie in a smaller disk and the supports of the shifted mollifiers
are contained in the original disk. The right-hand side converges locally
uniformly to $D_{z}^{\alpha }D_{w}^{\beta }M(z,w)$. Taking $\beta =\alpha $
and $w=z$, and expanding the squared norm of $D^{\alpha }X_{\varepsilon
}(z)-D^{\alpha }X_{\delta }(z)$, shows that the mollified maps and all their
derivatives are locally uniformly Cauchy. Since  $X_{\varepsilon
}\rightarrow X$ in $L_{\mathrm{loc}}^{2}$and equals $X$ almost everywhere,
the limit is a smooth representative of $X$, which we continue to denote by $%
X$. By continuity, the smooth representative is positive and has rank at
most one everywhere, since $\Vert X(z)\Vert _{\mathrm{HS}}=\rho _{1}(z)>0$.
Since $X(z_{0})$ is positive and nonzero, choose a fixed vector $e\in \ell
^{2}$ such that $\langle X(z_{0})e,e\rangle _{\ell ^{2}}>0$. By continuity,
there is a disk centered at $z_{0}$, contained in the preceding disk, where $%
\langle X(z)e,e\rangle _{\ell ^{2}}$ is positive. Define 
\begin{equation}
\widetilde{v}(z)=\frac{X(z)e}{\sqrt{\langle X(z)e,e\rangle _{\ell ^{2}}}}.
\label{eq:smooth-vectors}
\end{equation}%
The denominator is smooth and positive, so $\widetilde{v}$ is smooth. For a
positive operator of rank one, this formula gives $X(z)=\widetilde{v}(z)%
\widetilde{v}(z)^{\ast }$. Since we also have $X(z)=v(z)v^{\ast }(z)$ almost
everywhere, $v(z)$ and $\widetilde{v}(z)$ differ almost everywhere by a
scalar of modulus one. Consequently, replacing $\mathcal{K}(z,w)$ by $%
\langle \widetilde{v}(z),\widetilde{v}(w)\rangle _{\ell ^{2}}$, multiplies
each row and the corresponding column of every correlation determinant by
reciprocal phase factors, leaving it unchanged almost everywhere. We have
obtained the smooth kernel $\widetilde{\mathcal{K}}(z,w)=\langle \widetilde{v%
}(z),\widetilde{v}(w)\rangle _{\ell ^{2}}$. 

\medskip \refstepcounter{proofstep}\label{step:factor} \noindent \textbf{%
Step \theproofstep.} The determinants formed from the smooth correlation
kernel obtained in Step~\ref{step:smooth}, $\widetilde{\mathcal{K}}%
(z,w)=\langle \widetilde{v}(z),\widetilde{v}(w)\rangle _{\ell ^{2}}$, agree
pointwise with the continuous two- and three-point intensities of Lemma~\ref%
{lem:nearby-intensities} at nearby distinct pairs and triples. Indeed, they
agree almost everywhere and both sides are continuous. We can therefore
Taylor-expand these determinants and compare their coefficients with %
\eqref{eq:two-nearby} and \eqref{eq:three-nearby}. From now on, write $v$
for $\widetilde{v}$ and $\mathcal{K}$ for $\widetilde{\mathcal{K}}$. We now
prove the local factorization \eqref{eq:kernel-factor}. At a fixed point $z$%
, let $P^{\perp }$ be the orthogonal projection onto the subspace
perpendicular to $v(z)$, and set 
\begin{equation*}
u(z)=P^{\perp }\partial v(z),\qquad w(z)=P^{\perp }\overline{\partial }v(z).
\end{equation*}%
The two-point intensity is the Gram determinant of $v(z)$ and $v(z+\delta )$%
. Removing the component of the second vector along $v(z)$ gives 
\begin{equation*}
\rho _{2}(z,z+\delta )=\Vert v(z)\Vert ^{2}\Vert P^{\perp }v(z+\delta )\Vert
^{2}.
\end{equation*}%
Since $P^{\perp }v(z)=0$, Taylor's formula gives 
\begin{eqnarray*}
\rho _{2}(z,z+\delta ) &=&\Vert v(z)\Vert ^{2}\Vert u(z)\delta +w(z)%
\overline{\delta }\Vert ^{2}+O(|\delta |^{3}) \\
&=&\Vert v(z)\Vert ^{2}\bigl((\Vert u(z)\Vert ^{2}+\Vert w(z)\Vert
^{2})|\delta |^{2}+2\mathrm{Re}(\langle u(z),w(z)\rangle \delta ^{2})\bigr)%
+O(|\delta |^{3}).
\end{eqnarray*}%
The term involving $\delta ^{2}$ vanishes because the coefficient in %
\eqref{eq:two-nearby} is strictly positive and independent of the direction
of $\delta $. Thus, from the $2$-point identities, we have obtained 
\begin{equation}
\langle u(z),w(z)\rangle =0,\qquad \Vert v(z)\Vert ^{2}(\Vert u(z)\Vert
^{2}+\Vert w(z)\Vert ^{2})>0.  \label{eq:orthogonality}
\end{equation}

To factor a holomorphic factor out of $v(z)$ we now use the $3$-point
identities. For real $\epsilon $, set 
\begin{equation*}
r_{1}=P^{\perp }v(z+\epsilon ),\qquad r_{2}=P^{\perp }v(z+i\epsilon ).
\end{equation*}%
The three-point intensity is the Gram determinant of $v(z),v(z+\epsilon
),v(z+i\epsilon )$. Removing the components of the last two vectors along $%
v(z)$ leaves the determinant unchanged. The remaining vectors $r_{1},r_{2}$
are perpendicular to $v(z)$, so the Gram matrix becomes block diagonal: 
\begin{equation*}
\rho _{3}(z,z+\epsilon ,z+i\epsilon )=\Vert v(z)\Vert ^{2}\bigl(\Vert
r_{1}\Vert ^{2}\Vert r_{2}\Vert ^{2}-|\langle r_{1},r_{2}\rangle |^{2}\bigr).
\end{equation*}%
Taylor's formula gives 
\begin{equation*}
r_{1}=\epsilon (u+w)+O(\epsilon ^{2}),\qquad r_{2}=i\epsilon
(u-w)+O(\epsilon ^{2}).
\end{equation*}%
Substituting these expressions into the determinant gives 
\begin{equation*}
\rho _{3}(z,z+\epsilon ,z+i\epsilon )=4\epsilon ^{4}\Vert v(z)\Vert ^{2}%
\bigl(\Vert u(z)\Vert ^{2}\Vert w(z)\Vert ^{2}-|\langle u(z),w(z)\rangle
|^{2}\bigr)+O(|\epsilon |^{5}).
\end{equation*}%
On the other hand, \eqref{eq:three-nearby} bounds this intensity by $%
O(|\epsilon |^{6})$. Divide by $\epsilon ^{4}$ and let $\epsilon \rightarrow
0$. Using \eqref{eq:orthogonality}, we obtain 
\begin{equation*}
\Vert u(z)\Vert ^{2}\Vert w(z)\Vert ^{2}=0.
\end{equation*}%
Thus exactly one of $u(z),w(z)$ vanishes at each point. The sets where $%
\Vert u(z)\Vert >0$ and where $\Vert w(z)\Vert >0$ are disjoint open sets
whose union is the disk. Connectedness shows that the same choice holds
throughout the disk. Suppose $w=0$. Then $\overline{\partial }v(z)=b(z)v(z)$%
, where 
\begin{equation*}
b(z)=\frac{\langle \overline{\partial }v(z),v(z)\rangle _{\ell ^{2}}}{\Vert
v(z)\Vert _{\ell ^{2}}^{2}}
\end{equation*}%
is smooth. On a smaller disk the Cauchy--Green formula gives a smooth
solution $B(z)$ of $\overline{\partial }B=b(z)$. Let $f(z)=e^{-B(z)}v(z)$
and $a(z)=e^{B(z)}$. Direct differentiation gives 
\begin{equation*}
\overline{\partial }f(z)=e^{-B(z)}\bigl(\overline{\partial }v(z)-(b(z))v(z)%
\bigr)=0.
\end{equation*}%
Hence $f(z)$ is holomorphic and $a(z)$ is smooth and never zero. Taking 
\begin{equation*}
H(z,w)=\langle f(z),f(w)\rangle _{\ell ^{2}}=\frac{\langle v(z),v(w)\rangle
_{\ell ^{2}}}{a(z)\overline{a(w)}}
\end{equation*}%
gives a positive Gram kernel, holomorphic in the first variable and
antiholomorphic in the second, and the smooth kernel obtained in the
previous step satisfies $\mathcal{K}(z,w)=a(z)\overline{a(w)}H(z,w)$. If $u=0
$, simply replace $\mathcal{K}$ by its complex conjugate. This preserves all
correlation determinants and interchanges the roles of $u$ and $w$. We have
completed the proof of \eqref{eq:kernel-factor}. \hfill $\square $

\subsection{The curvature Lemma}

\label{sec:curvature}

\begin{lemma}
\label{lem:curvature} Let $F$ belong to the Peres--Vir\'{a}g class, and let $%
\rho _{1}(z)$ be its first zero intensity relative to $d\mu $. Then 
\begin{equation}
\partial \overline{\partial }\log \rho _{1}(z)=2\rho _{1}(z),\qquad z\in
\Omega .  \label{eq:liouville}
\end{equation}%
Thus the conformal metric $ds^{2}=\rho _{1}(z)|dz|^{2}$ has constant
Gaussian curvature $-4$.
\end{lemma}

\textbf{Proof. }

\setcounter{proofstep}{0} \medskip \refstepcounter{proofstep}\label%
{step:quadratic} \noindent\textbf{Step \theproofstep.} In this step we
establish the compatibility identity, valid near the diagonal, 
\begin{equation}  \label{eq:compatibility}
-\partial_z\overline{\partial}_z\log(1-G(z,w))=p(z),\qquad p(z)=\frac{%
(2\rho_1(z)+\partial\overline{\partial}\log \rho_1(z))^2}{8\rho_1(z)}>0.
\end{equation}
We start by observing that, when computing the normalized pair correlation
function $G(z,w)$ with the kernel \eqref{eq:kernel-factor} obtained in Lemma~%
\ref{lem:local-factorization}, the scalar factors cancel, giving: 
\begin{equation}  \label{eq:pair-ratio}
G(z,w)=1-\frac{|H(z,w)|^2}{H(z,z)H(w,w)}.
\end{equation}
Since $H(z,z)>0$, we may shrink the disk so that $H(z,w)\ne0$ for every pair
in it. For fixed $w$, $\log|H(z,w)|^2$ is harmonic in $z$. Taking logarithms
of $1-G(z,w)$ in \eqref{eq:pair-ratio} and then differentiating gives 
\begin{equation*}
-\partial_z\overline{\partial}_z\log(1-G(z,w)) =\partial\overline{\partial}%
\log H(z,z).
\end{equation*}
To identify the right-hand side, expand the logarithm of %
\eqref{eq:pair-ratio} at $w=z+\delta$. Hermitian symmetry of $H$ and its
holomorphic--antiholomorphic dependence cancel the constant, linear and pure
quadratic terms. The mixed term gives 
\begin{equation*}
-\log(1-G(z,z+\delta)) =\bigl(\partial\overline{\partial}\log H(z,z)\bigr)%
|\delta|^2 +O(|\delta|^3).
\end{equation*}
Since $G(z,z+\delta)=O(|\delta|^2)$, the quadratic coefficient of $%
-\log(1-G) $ is also that of $G$. Dividing \eqref{eq:two-nearby} by $%
\rho_1(z)\rho_1(z+\delta)$ gives 
\begin{equation*}
G(z,z+\delta) =\frac{(2\rho_1(z)+\partial\overline{\partial}\log\rho_1(z))^2%
} {8\rho_1(z)}|\delta|^2+O(|\delta|^3).
\end{equation*}
Comparison identifies $\partial\overline{\partial}\log H(z,z)$ with $p(z)$
and proves \eqref{eq:compatibility}. Its positivity follows from %
\eqref{eq:second-conditional}.

\medskip \refstepcounter{proofstep}\label{step:quartic} \noindent\textbf{%
Step \theproofstep.} In this step we prove \eqref{eq:liouville} throughout $%
\Omega$ by first computing the coefficient of $z^2\overline{z}^2$ in 
\begin{equation}  \label{eq:normalized-pair-at-zero}
G(0,z)=\frac{\rho_2(0,z)}{\rho_1(0)\rho_1(z)},
\end{equation}
using the Kac--Rice formula \eqref{eq:kac-rice} for the two-point intensity $%
\rho_2$, and then comparing it with the coefficient of $z^2\overline{z}^2$
obtained from the case $w=0$ of the differential identity %
\eqref{eq:compatibility}.

Fix a point and translate it to zero. Since $K(z,0)$ is nonzero nearby,
replace $F$ by 
\begin{equation}
z\longmapsto \frac{\sqrt{K(0,0)}F(z)}{K(z,0)}.  \label{replacement}
\end{equation}%
The new covariance satisfies $K(z,0)=1$. If $c=\partial _{z}\partial _{%
\overline{w}}K(0,0)>0$, use the local coordinate 
\begin{equation*}
\zeta (z)=c^{-1/2}\left. \partial _{\overline{w}}K(z,w)\right\vert _{w=0}.
\end{equation*}%
Its derivative at zero is $\sqrt{c}>0$, so it has a holomorphic local
inverse. The chain rule then gives the following identities for the
covariance in the new coordinate, again denoted by $z$: 
\begin{equation}
K(z,0)=1,\qquad \left. \partial _{\overline{w}}K(z,w)\right\vert _{w=0}=z.
\label{eq:exact-normal}
\end{equation}%
Taking complex conjugates also gives $\partial _{z}K(0,w)=\overline{w}$.
These identities determine every coefficient for which one exponent is zero
or one. The expansion through total degree six is therefore 
\begin{eqnarray}
K(z,w) &=&1+z\overline{w}+\beta z^{2}\overline{w}^{2}+\eta z^{3}\overline{w}%
^{2}+\overline{\eta }z^{2}\overline{w}^{3}  \notag \\
&&+\gamma z^{3}\overline{w}^{3}+\chi z^{4}\overline{w}^{2}+\overline{\chi }%
z^{2}\overline{w}^{4}+O_{\mathrm{tot}}(7),  \label{eq:normal-series}
\end{eqnarray}%
with $\gamma $ real and $\beta >0$. In fact, the constant and linear
coefficients of $F$ are independent standard complex Gaussians, orthogonal
to every higher coefficient, and $\beta $ is the variance of the quadratic
coefficient. Its positivity follows from \eqref{eq:derivative-assumption}.
The notation $O_{\mathrm{tot}}(m)$ means terms of total degree at least $m$.
Thus (\ref{eq:normal-series}) assures that expanding $K$ through total
degree six is sufficient to compute the coefficients $A_{4}$ and $C_{2}$
below.

The replacement (\ref{replacement}) preserves the zeros of $F(z)$, since the
multiplier is holomorphic and never zero on the chosen neighborhood. The
ratio $G $ is unchanged by the coordinate change, and the kernel still has
the form \eqref{eq:kernel-factor}. We may make the comparison in these
coordinates. Write $t=|z|^{2}$ and 
\begin{eqnarray}
K(z,z) &=&1+tU(z),  \notag \\
U(z) &=&1+\beta t+\eta z^{2}\overline{z}+\overline{\eta }z\overline{z}%
^{2}+\gamma t^{2}+\chi z^{3}\overline{z}+\overline{\chi }z\overline{z}%
^{3}+O_{\mathrm{tot}}(5).  \label{eq:U}
\end{eqnarray}%
For the pair $(0,z)$, the three matrices in Kac--Rice are 
\begin{equation*}
A=\left( 
\begin{array}{cc}
1 & 1 \\ 
1 & K(z,z)%
\end{array}%
\right) ,\qquad B=\left( 
\begin{array}{cc}
0 & \overline{z} \\ 
0 & \partial K(z,z)%
\end{array}%
\right) ,\qquad C=\left( 
\begin{array}{cc}
1 & 1 \\ 
1 & \partial \overline{\partial }K(z,z)%
\end{array}%
\right) .
\end{equation*}%
Since $\det A=tU(z)$, computing the conditional covariance matrix %
\eqref{eq:conditional-covariance} gives 
\begin{equation}
\Lambda =\left( 
\begin{array}{cc}
(U(z)-1)/U(z) & -\overline{z}\,\overline{\partial }\log U(z) \\ 
-z\,\partial \log U(z) & tU(z)\,\partial \overline{\partial }\log U(z)%
\end{array}%
\right) .  \label{eq:conditional-matrix}
\end{equation}%
It follows that 
\begin{equation}
\rho _{2}(0,z)=\frac{(U(z)-1)\partial \overline{\partial }\log
U(z)+|\partial \log U(z)|^{2}}{U(z)}.  \label{eq:exact-pair-normal}
\end{equation}

To compute the required coefficients, expand 
\begin{equation*}
\log U(z)=\beta t+\eta z^{2}\overline{z}+\overline{\eta }z\overline{z}%
^{2}+(\gamma -\beta ^{2}/2)t^{2}+\chi z^{3}\overline{z}+\overline{\chi }z%
\overline{z}^{3}+O_{\mathrm{tot}}(5).
\end{equation*}%
The coefficients of $t^{2}$ in $(U(z)-1)\partial \overline{\partial }\log
U(z)$ and $|\partial \log U(z)|^{2}$ are, respectively, 
\begin{equation*}
5\beta \gamma -2\beta ^{3}+4|\eta |^{2},\qquad 4\beta \gamma -2\beta
^{3}+5|\eta |^{2}.
\end{equation*}%
The leading term in \eqref{eq:exact-pair-normal} is $2\beta ^{2}t$. Dividing
the sum of the two fourth-order coefficients by $U(z)=1+\beta t+O(|z|^{3})$
therefore subtracts a further $2\beta ^{3}$. Thus the coefficient of $t^{2}$
in $\rho _{2}(0,z)$ is $9\beta \gamma +9|\eta |^{2}-6\beta ^{3}$. The
Edelman--Kostlan formula applied to $K(z,z)=1+tU(z)$ gives 
\begin{eqnarray}
\rho _{1}(z) &=&1+(4\beta -2)t+6\eta z^{2}\overline{z}+6\overline{\eta }z%
\overline{z}^{2}  \notag \\
&&+9(\gamma -\beta +1/3)t^{2}+8\chi z^{3}\overline{z}+8\overline{\chi }z%
\overline{z}^{3}+O_{\mathrm{tot}}(5).  \label{eq:g-expansion}
\end{eqnarray}%
In particular, $\rho _{1}(0)=1$. This allows us to conclude, from %
\eqref{eq:normalized-pair-at-zero} and \eqref{eq:exact-pair-normal}, that 
\begin{equation*}
G(0,z)=\frac{\rho _{2}(0,z)}{\rho _{1}(z)}=2\beta ^{2}|z|^{2}+O(|z|^{3}).
\end{equation*}%
To find its fourth-order coefficient, divide by $\rho _{1}(z)=1+(4\beta
-2)t+O(|z|^{3})$. This subtracts $2\beta ^{2}(4\beta -2)$ from the
coefficient of $t^{2}$ in $\rho _{2}$. Consequently the coefficient of $z^{2}%
\overline{z}^{2}$ in $G(0,z)$ is 
\begin{equation}
A_{4}=9\beta \gamma +9|\eta |^{2}+4\beta ^{2}-14\beta ^{3}.
\label{eq:actual-quartic}
\end{equation}%
We next compute the coefficient of $z\overline{z}$ in $p(z)$. Taking the
logarithm in \eqref{eq:g-expansion} and differentiating gives 
\begin{eqnarray*}
2\rho _{1}(z)+\partial \overline{\partial }\log \rho _{1}(z) &=&4\beta
+12\eta z+12\overline{\eta }\,\overline{z} \\
&&+(36\gamma +4\beta -32\beta ^{2})t+24\chi z^{2}+24\overline{\chi }\,%
\overline{z}^{2}+O_{\mathrm{tot}}(3).
\end{eqnarray*}%
The coefficient of $t$ in the square of this expression is $8\beta (36\gamma
+4\beta -32\beta ^{2})+288|\eta |^{2}$. First divide this coefficient by $8$%
. The factor $1/\rho _{1}(z)$ then subtracts $2\beta ^{2}(4\beta -2)$.
Therefore the coefficient of $z\overline{z}$ in $p(z)$ is 
\begin{equation}
C_{2}=36\beta \gamma +36|\eta |^{2}+8\beta ^{2}-40\beta ^{3}.
\label{eq:derivatives-p}
\end{equation}

Since $G$ is symmetric, \eqref{eq:compatibility} with $w=0$ applies to $%
G(0,z)$ as well. Now 
\begin{equation*}
-\log (1-G(0,z))=G(0,z)+\frac{1}{2}G(0,z)^{2}+O(|z|^{6}).
\end{equation*}%
Its coefficient of $z^{2}\overline{z}^{2}$ is $A_{4}+2\beta ^{4}$. Applying $%
\partial \overline{\partial }$ multiplies this coefficient by four, so %
\eqref{eq:compatibility} implies 
\begin{equation*}
C_{2}=4(A_{4}+2\beta ^{4}).
\end{equation*}%
Substituting \eqref{eq:actual-quartic}--\eqref{eq:derivatives-p} yields 
\begin{equation}
0=4A_{4}+8\beta ^{4}-C_{2}=8\beta ^{2}(\beta -1)^{2}.
\label{eq:quartic-obstruction}
\end{equation}%
Nondegeneracy gives $\beta >0$, so $\beta =1$. The terms containing $\gamma $
and $|\eta |^{2}$ cancel. The terms containing $\chi $ do not enter the
indicated coefficients.

At the origin, \eqref{eq:g-expansion} gives $(\partial\overline{\partial}%
\log \rho_1)(0)=4\beta-2$ and $\rho_1(0)=1$. It follows that $%
(\partial\overline\partial\log \rho_1)(0)=2\rho_1(0)$. The point was
arbitrary. Under a holomorphic change of variable $\psi$, 
\begin{eqnarray*}
\rho_{1,\psi}(z) &=& |\psi^{\prime }(z)|^2\rho_1(\psi(z)), \\
\partial\overline{\partial}\log \rho_{1,\psi}(z) &=& |\psi^{\prime
}(z)|^2(\partial\overline{\partial}\log \rho_1)(\psi(z)),
\end{eqnarray*}
because $\log|\psi^{\prime }(z)|^2$ is harmonic. The equality therefore
holds in the original coordinate, proving \eqref{eq:liouville} throughout $%
\Omega$.

\hfill$\square$

\subsection{Proof of Theorem~\protect\ref{thm:main}}

\label{sec:main-proof}

The converse follows from Peres--Vir\'{a}g~\cite{PV}. Indeed, a GAF with
covariance \eqref{eq:main-covariance} has the same law as $h(F_{\mathbb{D}%
}\circ \phi )$. Multiplication by $h$ does not change its zeros. Restricting
the Peres--Vir\'{a}g process to $\phi (\Omega )$ and changing variables by $%
\phi $ therefore gives the DPP kernel \eqref{eq:main-dpp-kernel}. To verify %
\eqref{eq:derivative-assumption}, consider the independent Gaussian
coefficients multiplying $h(z),h(z)\phi (z),h(z)\phi ^{2}(z)$. At every
point $z$ their values and first two derivatives form a matrix with
determinant 
\begin{equation*}
\det \bigl[(h\phi ^{j})^{(i)}(z)\bigr]_{i,j=0}^{2}=2h(z)^{3}\phi ^{\prime
}(z)^{3}\neq 0.
\end{equation*}%
These three independent coefficients give a positive definite contribution
to the covariance matrix in \eqref{eq:derivative-assumption}. The remaining
coefficients add a positive semidefinite matrix, so the full covariance is
positive definite.

We now prove necessity. Assume that $F$ belongs to the Peres--Vir\'ag class.
By Lemma~\ref{lem:curvature}, its first intensity satisfies %
\eqref{eq:liouville}. We first obtain the global covariance representation
and then prove that its map $\phi$ is one-to-one.

\setcounter{proofstep}{0} \medskip \refstepcounter{proofstep}\label%
{step:global} \noindent \textbf{Step \theproofstep.} We first obtain the
covariance kernel on a sufficiently small disk $U\subseteq \Omega $.
Liouville's theorem, Lemma~\ref{lem:liouville}, applied with $g=\rho _{1}$
to \eqref{eq:liouville}, gives a holomorphic map $\phi _{0}:U\rightarrow 
\mathbb{D}$ with $\phi _{0}^{\prime }\neq 0$ and 
\begin{equation*}
\rho _{1}(z)=\frac{|\phi _{0}^{\prime }(z)|^{2}}{(1-|\phi _{0}(z)|^{2})^{2}}.
\end{equation*}%
The right-hand side is the first zero intensity of $F_{\mathbb{D}}\circ \phi
_{0}$, by the Edelman--Kostlan formula. Thus this GAF and $F|_{U}$ have the
same first zero intensity. Sodin's Calabi rigidity theorem~\cite[Theorem 2]%
{Sodin}, in the form stated in~\cite[Theorem 2.5.2, pp.~30--31]{HKPV}, gives
a zero-free holomorphic function $h_{0}$ on $U$ such that 
\begin{equation}
K(z,w)=\frac{h_{0}(z)\overline{h_{0}(w)}}{1-\phi _{0}(z)\overline{\phi
_{0}(w)}},\qquad z,w\in U.  \label{eq:local-covariance}
\end{equation}

It remains to extend this representation to the whole domain $\Omega $. Fix $%
z_{0}\in U$. By M\"{o}bius invariance of the disk, we may assume $\phi
_{0}(z_{0})=0$, after a corresponding adjustment of $h_{0}$ in %
\eqref{eq:local-covariance}. Multiplying $h_{0}$ by a constant of modulus
one, we may also assume 
\begin{equation*}
h_{0}(z_{0})=\sqrt{K(z_{0},z_{0})}>0.
\end{equation*}%
Evaluating \eqref{eq:local-covariance} at $w=z_{0}$ now expresses $h_{0}(z)$
as a constant multiple of the globally defined holomorphic function $%
K(z,z_{0})$. Hence $h_{0}$ extends holomorphically to $\Omega $. Now
differentiate \eqref{eq:local-covariance} in $\overline{w}$ and set $w=z_{0}$%
. Since $\phi _{0}(z_{0})=0$, this identity expresses $h_{0}(z)\phi _{0}(z)$
as a constant linear combination of $K(z,z_{0})$ and $\left. \partial _{%
\overline{w}}K(z,w)\right\vert _{w=z_{0}}$. The coefficient by which we
divide is $h_{0}(z_{0})\overline{\phi _{0}^{\prime }(z_{0})}\neq 0$. Thus $%
h_{0}\phi _{0}$ also extends holomorphically to $\Omega $. Denote the two
extensions by $h$ and $j$, respectively. The diagonal of %
\eqref{eq:local-covariance} gives, on $U$, the identity 
\begin{equation}
K(z,z)\bigl(|h(z)|^{2}-|j(z)|^{2}\bigr)=|h(z)|^{4}.
\label{eq:continued-diagonal-identity}
\end{equation}%
Both sides are real-analytic on the connected domain $\Omega $, including
points where $h$ vanishes. Since their difference vanishes on $U$,
real-analytic continuation proves \eqref{eq:continued-diagonal-identity}
throughout $\Omega $. Since $h(z_{0})>0$, the zeros of $h$ are isolated. At
a point where $h(z)\neq 0$, equation~\eqref{eq:continued-diagonal-identity}
and $K(z,z)>0$ imply 
\begin{equation*}
\left\vert \frac{j(z)}{h(z)}\right\vert <1.
\end{equation*}%
The quotient $j(z)/h(z)$ is therefore bounded near every zero of $h(z)$. The
removable singularity theorem extends it to a holomorphic function $\phi (z)$
on $\Omega $, which is nonconstant on $U$, where $\phi (z)=\phi _{0}(z)$ and 
$\phi _{0}^{\prime }(z)\neq 0$. By continuity $|\phi (z)|\leq 1$ everywhere
and, by the maximum modulus principle, $|\phi (z)|<1$ throughout $\Omega $.
On $\{z:h(z)\neq 0\}$, equation~\eqref{eq:continued-diagonal-identity} can
now be written as 
\begin{equation*}
|h(z)|^{2}=K(z,z)(1-|\phi (z)|^{2}).
\end{equation*}%
Both sides are continuous, so this equality extends across the isolated
zeros of $h$. Its right-hand side is strictly positive. Consequently $h$ has
no zeros. Taking logarithms gives 
\begin{equation*}
\log K(z,z)=\log |h(z)|^{2}-\log (1-|\phi (z)|^{2}).
\end{equation*}%
The first term is harmonic because $h$ is holomorphic and zero-free.
Applying the Edelman--Kostlan formula therefore gives 
\begin{equation*}
\rho _{1}(z)=\frac{|\phi ^{\prime }(z)|^{2}}{(1-|\phi (z)|^{2})^{2}}.
\end{equation*}%
Since $\rho _{1}(z)>0$, the derivative $\phi ^{\prime }(z)$ never vanishes.
The kernel $K(z,w)$ and the right-hand side of \eqref{eq:local-covariance},
with the extensions $h,\phi $, are now holomorphic in $z$ and
antiholomorphic in $w$ on $\Omega \times \Omega $. For each fixed $w\in U$,
the identity theorem in $z$ extends their equality from $U$ to $\Omega $.
For each fixed $z\in \Omega $, the same theorem in the second variable then
extends the equality from $U$ to $\Omega $. Thus \eqref{eq:local-covariance}
holds on the whole of $\Omega \times \Omega $, proving %
\eqref{eq:main-covariance}.

\medskip \refstepcounter{proofstep}\label{step:injective} \noindent \textbf{%
Step \theproofstep.} We show that the map $\phi $ obtained in Step~\ref%
{step:global} is one-to-one, arguing by contradiction. Suppose $\phi
(z_{1})=\phi (z_{2})$ with $z_{1}\neq z_{2}$. Since $\phi ^{\prime }(z)$
never vanishes, the inverse function theorem gives disjoint neighborhoods on
which the two inverse branches are defined. Restrict their images to the
same sufficiently small disk. We obtain disjoint relatively compact
neighborhoods $U,V$ such that $\phi |_{U}$ and $\phi |_{V}$ are one-to-one
and have that disk as their common image. Set 
\begin{equation*}
T=(\phi |_{V})^{-1}\circ (\phi |_{U}):U\longrightarrow V.
\end{equation*}%
Then $T$ is a holomorphic bijection $T:U\rightarrow V$ and $\phi (Tz)=\phi
(z)$. By \eqref{eq:main-covariance}, for every fixed $z\in U$, 
\begin{equation*}
\frac{F(z)}{h(z)}=\frac{F(Tz)}{h(Tz)}\qquad \text{almost surely}.
\end{equation*}%
Indeed, the covariance identity \eqref{eq:main-covariance} gives, for $z\in U
$, 
\begin{equation*}
\mathbb{E}\left\vert \frac{F(z)}{h(z)}-\frac{F(Tz)}{h(Tz)}\right\vert ^{2}=%
\frac{2}{1-|\phi (z)|^{2}}-2\mathrm{Re}\frac{1}{1-\phi (z)\overline{\phi (Tz)%
}}=0\text{.}
\end{equation*}%
Apply the almost sure equality on a countable dense subset of $U$ and
intersect the resulting events of probability one. By continuity, the
equality holds simultaneously for every $z\in U$ on the intersection event.
Since $h(z)$ never vanishes and $T=(\phi |_{V})^{-1}\circ (\phi |_{U})$ is a
biholomorphic map from $U$ onto $V$, the zeros in $U$ and $V$ correspond
bijectively. Let $N_{F}(E)$ denote the number of zeros in $E$. Consequently, 
\begin{equation}
N_{F}(U)=N_{F}(V)\qquad \text{almost surely}.  \label{almostSure}
\end{equation}%
By shrinking the common image disk and the corresponding neighborhoods $U,V$%
, we may assume that 
\begin{equation*}
m:=\mathbb{E}[N_{F}(U)]=\mathbb{E}[N_{F}(V)]=\int_{U}\rho _{1}(z)\,d\mu
(z),\qquad 0<m<1.
\end{equation*}%
Indeed, $\rho _{1}$ is continuous and strictly positive, and the
neighborhoods shrink to $z_{1},z_{2}$ as the radius of the common image disk
tends to zero. Since $U$ and $V$ are disjoint, integration of the (global
repulsive) two-point determinantal identity \eqref{eq:modulus-identity},
gives 
\begin{eqnarray*}
\mathbb{E}[N_{F}(U)N_{F}(V)] &=&\int_{U\times V}\rho _{2}(z,w)\,d\mu
(z)\,d\mu (w) \\
&\leq &\left( \int_{U}\rho _{1}(z)\,d\mu (z)\right) \left( \int_{V}\rho
_{1}(w)\,d\mu (w)\right) =m^{2}.
\end{eqnarray*}%
On the other hand, since $N_{F}(U)$ is a nonnegative integer, $%
N_{F}(U)^{2}\geq N_{F}(U)$. Together with the almost sure equality (\ref%
{almostSure}), this gives 
\begin{equation*}
\mathbb{E}[N_{F}(U)N_{F}(V)]=\mathbb{E}[N_{F}(U)^{2}]\geq \mathbb{E}%
[N_{F}(U)]=m>m^{2}.
\end{equation*}%
This contradiction proves that $\phi $ is one-to-one and completes the proof
of Theorem~\ref{thm:main}.

\section{Examples and comments}

\label{sec:examples}

\subsection{Restrictions and changes of variable}

Every one-to-one holomorphic map $\phi:\Omega\to\mathbb{D}$ gives an example
through $F_{\mathbb{D}}\circ\phi$. Taking $\phi(z)=z$ on a strict subdomain
of the disk shows why $\phi(\Omega)=\mathbb{D}$ is not required. In
particular, we may restrict the disk process to an annulus. The domain in
the theorem need not be simply connected.

\subsection{A function with one random zero}

\label{subsec:linear-example}

Let $F(z)=\xi _{0}+\xi _{1}z$ on $\mathbb{C}$. The zeros of $F(z)$ are
determinantal, with kernel 
\begin{equation*}
\mathcal{K}(z,w)=\frac{1}{(1+|z|^{2})(1+|w|^{2})}.
\end{equation*}%
The corresponding operator projects onto the one-dimensional space spanned
by $(1+|z|^{2})^{-1}$. Here $K(z,w)=1+z\overline{w}$ and $\partial \overline{%
\partial }\log \rho _{1}(z)=-2\rho _{1}(z)$. The two-point intensity
vanishes identically. This example is outside the Peres--Vir\'{a}g class,
since its covariance matrix is 
\begin{equation*}
\mathrm{Cov}\bigl(F(z),F^{\prime }(z),F^{\prime \prime }(z)\bigr)=\left( 
\begin{array}{ccc}
1+|z|^{2} & z & 0 \\ 
\overline{z} & 1 & 0 \\ 
0 & 0 & 0%
\end{array}%
\right) .
\end{equation*}%
This matrix has rank two and determinant zero at every point. It is
therefore not positive definite, so \eqref{eq:derivative-assumption} fails.

\subsection{Perturbing the constant coefficient of the Hardy GAF}

Consider the disk limit of the weighted Szeg\H{o} family studied by Katori
and Shirai~\cite{KS}, 
\begin{equation}
F_{r}(z)=\frac{\xi _{0}}{\sqrt{1+r}}+\sum_{n=1}^{\infty }\xi
_{n}z^{n},\qquad r>0.  \label{eq:perturbed-gaf}
\end{equation}%
One can directly check that the covariance of $(F_{r},F_{r}^{\prime
},F_{r}^{\prime \prime })$ is positive definite. At the origin, 
\begin{equation*}
(\partial \overline{\partial }\log \rho _{1,r})(0)-2\rho _{1,r}(0)=-4r\neq 0.
\end{equation*}%
Thus $\rho _{1,r}$ does not satisfy the necessary equation %
\eqref{eq:liouville} established in Lemma~\ref{lem:curvature}. Its zero
process is therefore not a Hermitian DPP for any $r>0$. \cite[Proposition 1.7%
]{KS} gives a nonempty interval of positive $r$ for which zeros are
negatively correlated. This shows that perturbing the constant coefficient
can preserve repulsion, but destroys the Hermitian determinantal property.

\subsection{Constant curvature without global injectivity}

\label{subsec:noninjective-example}

We give an example satisfying \eqref{eq:derivative-assumption} whose zeros
are determinantal on sufficiently small disks, but fail to be determinantal
globally. Let $\Omega =\{z\in \mathbb{C}:\Re z<0\}$, $\phi (z)=e^{z}$ and
define 
\begin{equation}
F(z)=F_{\mathbb{D}}(e^{z})=\sum_{n=0}^{\infty }\xi _{n}e^{nz},\qquad z\in
\Omega .  \label{eq:exponential-gaf}
\end{equation}%
Since $|e^{z}|<1$ on $\Omega $, this is a centered circular GAF with
covariance 
\begin{equation}
K(z,w)=\frac{1}{1-e^{z+\overline{w}}},\qquad z,w\in \Omega .
\label{eq:exponential-covariance}
\end{equation}%
The map $\phi $ is locally one-to-one, since $\phi ^{\prime }(z)=e^{z}\neq 0$%
. The covariance of $(F(z),F^{\prime }(z),F^{\prime \prime }(z))$ is
positive definite at every point. Writing $x=\Re z$, direct differentiation
of the Edelman--Kostlan formula applied to \eqref{eq:exponential-covariance}
shows that 
\begin{equation}
\partial \overline{\partial }\log \rho _{1}(z)=2\rho _{1}(z).
\label{eq:exponential-intensity}
\end{equation}%
Thus the metric $ds^{2}=\rho _{1}(z)|dz|^{2}$ has Gaussian curvature $-4$
throughout $\Omega $. On every sufficiently small disk $U\subset \Omega $,
the exponential map is one-to-one. Restricting the Peres--Vir\'{a}g process
to $e^{U}$ and changing variables by $\phi $ therefore shows that the zeros
of $F|_{U}$ form a DPP, by~\cite[Theorem 1]{PV}. Its correlation kernel
relative to $d\mu $ is 
\begin{equation*}
\mathcal{K}_{U}(z,w)=\frac{e^{z+\overline{w}}}{(1-e^{z+\overline{w}})^{2}}%
,\qquad z,w\in U.
\end{equation*}%
However, on the whole domain, $\phi (z+2\pi i)=\phi (z)$, so $\phi $ is not
globally one-to-one. Correspondingly, \eqref{eq:exponential-gaf} gives 
\begin{equation}
F(z+2\pi i)=F(z),\qquad z\in \Omega .  \label{eq:exponential-periodicity}
\end{equation}%
By the argument of Step~\ref{step:injective} in the proof of Theorem~\ref%
{thm:main}, these forced pairs of zeros prevent a Hermitian determinantal
representation on $\Omega $.

\section*{Acknowledgement}

The author is thankful to OpenAI for the award of an `Academic Researchers
Account', providing free access to models aligned with scientific research.
We share information of potential methodological interest for mathematical
research. A substantial use of the AI tool was made during the experimental
part of this work, where a large database of kernels was tested, initially
with the goal of finding a new GAF with determinantal zeros. After several
frustrated attempts, `the curvature hypothesis' was tested on the same
database, which included constant curvature examples such as Bergman
kernels, as well as exotic $q$-analogues. The experimental results were then
organized in a knowledge graph, which was subsequently used in the
theoretical investigation, to retrieve connected examples, verify
computations and develop some technical arguments (which were then rederived
and simplified by the author), and also in the final revisions, to enhance
clarity of exposition. This research was funded by the Austrian Science Fund
(FWF) through the project 10.55776/PAT8205923.

\bigskip {\small \noindent Faculty of Mathematics, University of Vienna,%
\newline
Oskar-Morgenstern-Platz 1, A-1090, Vienna, Austria\newline
\textit{E-mail address:} \texttt{abreuluisdaniel@gmail.com} }

\end{document}